\documentclass[12pt,twoside,leqno]{article}
\usepackage{amsmath}
\usepackage{amssymb}
\usepackage{amsxtra}
\usepackage{amscd}
\usepackage{amsthm}
\usepackage[mathscr]{eucal}

\theoremstyle{plain}
\newtheorem{thm}[subsection]{Theorem}
\newtheorem{prop}[subsection]{Proposition}

\newtheorem{lem}[subsection]{Lemma}
\newtheorem{conj}[subsection]{Conjecture}

\theoremstyle{definition}

\newtheorem{para}[subsection]{}

\newenvironment{pf}{\proof[\proofname]}{\endproof}

\begin{document}

\title{Deligne--Beilinson cohomology and log Hodge theory}

\author{Kazuya Kato, Chikara Nakayama, Sampei Usui}

\maketitle

\newcommand\Cal{\mathcal}
\newcommand\define{\newcommand}

\define\gp{\mathrm{gp}}%
\define\fs{\mathrm{fs}}%
\define\an{\mathrm{an}}%
\define\mult{\mathrm{mult}}%
\define\Ker{\mathrm{Ker}\,}%
\define\Coker{\mathrm{Coker}\,}%
\define\Hom{\mathrm{Hom}\,}%
\define\Ext{\mathrm{Ext}\,}%
\define\rank{\mathrm{rank}\,}%
\define\gr{\mathrm{gr}}%
\define\cHom{\Cal{Hom}}
\define\cExt{\Cal Ext\,}%

\define\cA{\Cal A}
\define\cB{\Cal B}
\define\cC{\Cal C}
\define\cD{\Cal D}
\define\cO{\Cal O}
\define\cS{\Cal S}
\define\cM{\Cal M}
\define\cG{\Cal G}
\define\cH{\Cal H}
\define\cE{\Cal E}
\define\cF{\Cal F}
\define\cN{\Cal N}
\define\cR{\Cal R}
\define\cV{\Cal V}
\define\fF{\frak F}
\define\Dc{\check{D}}
\define\Ec{\check{E}}

\newcommand{\N}{{\mathbb{N}}}
\newcommand{\Q}{{\mathbb{Q}}}
\newcommand{\Z}{{\mathbb{Z}}}
\newcommand{\R}{{\mathbb{R}}}
\newcommand{\C}{{\mathbb{C}}}
\newcommand{\bN}{{\mathbb{N}}}
\newcommand{\bQ}{{\mathbb{Q}}}
\newcommand{\bF}{{\mathbb{F}}}
\newcommand{\bZ}{{\mathbb{Z}}}
\newcommand{\bP}{{\mathbb{P}}}
\newcommand{\bR}{{\mathbb{R}}}
\newcommand{\bC}{{\mathbb{C}}}
\newcommand{\bbQ}{{\bar \mathbb{Q}}}
\newcommand{\ol}[1]{\overline{#1}}
\newcommand{\too}{\longrightarrow}
\newcommand{\respect}{\rightsquigarrow}
\newcommand{\compatible}{\leftrightsquigarrow}
\newcommand{\upc}[1]{\overset {\lower 0.3ex \hbox{${\;}_{\circ}$}}{#1}}
\newcommand{\Gmlog}{\bG_{m, \log}}
\newcommand{\Gm}{\bG_m}
\newcommand{\ep}{\varepsilon}
\newcommand{\Spec}{\operatorname{Spec}}
\newcommand{\nilp}{\operatorname{nilp}}
\newcommand{\prim}{\operatorname{prim}}
\newcommand{\val}{{\mathrm{val}}} 
\newcommand{\n}{\operatorname{naive}}
\newcommand{\bs}{\operatorname{\backslash}}
\newcommand{\Gal}{\operatorname{{Gal}}}
\newcommand{\gal}{{\rm {Gal}}({\bar \Q}/{\Q})}
\newcommand{\galp}{{\rm {Gal}}({\bar \Q}_p/{\Q}_p)}
\newcommand{\gall}{{\rm{Gal}}({\bar \Q}_\ell/\Q_\ell)}
\newcommand{\wep}{W({\bar \Q}_p/\Q_p)}
\newcommand{\wel}{W({\bar \Q}_\ell/\Q_\ell)}
\newcommand{\Ad}{{\rm{Ad}}}
\newcommand{\BS}{{\rm {BS}}}
\newcommand{\even}{\operatorname{even}}
\newcommand{\End}{{\rm {End}}}
\newcommand{\odd}{\operatorname{odd}}
\newcommand{\GL}{\operatorname{GL}}
\newcommand{\np}{\text{non-$p$}}
\newcommand{\g}{{\gamma}}
\newcommand{\G}{{\Gamma}}
\newcommand{\Lam}{{\Lambda}}
\newcommand{\La}{{\Lambda}}
\newcommand{\lam}{{\lambda}}
\newcommand{\la}{{\lambda}}
\newcommand{\uL}{{{\hat {L}}^{\rm {ur}}}}
\newcommand{\uQp}{{{\hat \Q}_p}^{\text{ur}}}
\newcommand{\sel}{\operatorname{Sel}}
\newcommand{\dt}{{\rm{Det}}}
\newcommand{\Sig}{\Sigma}
\newcommand{\fil}{{\rm{fil}}}
\newcommand{\SL}{{\rm{SL}}}
\newcommand{\spl}{{\rm{spl}}}
\newcommand{\st}{{\rm{st}}}
\newcommand{\Isom}{{\rm {Isom}}}
\newcommand{\Mor}{{\rm {Mor}}}
\newcommand{\bg}{\bar{g}}
\newcommand{\id}{{\rm {id}}}
\newcommand{\cone}{{\rm {cone}}}
\newcommand{\al}{a}
\newcommand{\ChL}{{\cal{C}}(\La)}
\newcommand{\Image}{{\rm {Image}}}
\newcommand{\toric}{{\operatorname{toric}}}
\newcommand{\torus}{{\operatorname{torus}}}
\newcommand{\Aut}{{\rm {Aut}}}
\newcommand{\Qp}{{\mathbb{Q}}_p}
\newcommand{\barQp}{{\mathbb{Q}}_p}
\newcommand{\Qpur}{{\mathbb{Q}}_p^{\rm {ur}}}
\newcommand{\Zp}{{\mathbb{Z}}_p}
\newcommand{\Zl}{{\mathbb{Z}}_l}
\newcommand{\Ql}{{\mathbb{Q}}_l}
\newcommand{\Qlur}{{\mathbb{Q}}_l^{\rm {ur}}}
\newcommand{\F}{{\mathbb{F}}}
\newcommand{\eps}{{\epsilon}}
\newcommand{\epsLa}{{\epsilon}_{\La}}
\newcommand{\epsLaVxi}{{\epsilon}_{\La}(V, \xi)}
\newcommand{\epsOLaVxi}{{\epsilon}_{0,\La}(V, \xi)}
\newcommand{\Qplin}{{\mathbb{Q}}_p(\mu_{l^{\infty}})}
\newcommand{\otimesQplin}{\otimes_{\Qp}{\mathbb{Q}}_p(\mu_{l^{\infty}})}
\newcommand{\galFl}{{\rm{Gal}}({\bar {\Bbb F}}_\ell/{\Bbb F}_\ell)}
\newcommand{\gallur}{{\rm{Gal}}({\bar \Q}_\ell/\Q_\ell^{\rm {ur}})}
\newcommand{\galFF}{{\rm {Gal}}(F_{\infty}/F)}
\newcommand{\galFv}{{\rm {Gal}}(\bar{F}_v/F_v)}
\newcommand{\galF}{{\rm {Gal}}(\bar{F}/F)}
\newcommand{\epsVxi}{{\epsilon}(V, \xi)}
\newcommand{\epsOVxi}{{\epsilon}_0(V, \xi)}
\newcommand{\plim}{\lim_
{\scriptstyle 
\longleftarrow \atop \scriptstyle n}}
\newcommand{\sig}{{\sigma}}
\newcommand{\ga}{{\gamma}}
\newcommand{\del}{{\delta}}
\newcommand{\Vss}{V^{\rm {ss}}}
\newcommand{\Bst}{B_{\rm {st}}}
\newcommand{\Dpst}{D_{\rm {pst}}}
\newcommand{\Dcrys}{D_{\rm {crys}}}
\newcommand{\DdR}{D_{\rm {dR}}}
\newcommand{\Fin}{F_{\infty}}
\newcommand{\Kla}{K_{\lambda}}
\newcommand{\Ola}{O_{\lambda}}
\newcommand{\Mla}{M_{\lambda}}
\newcommand{\Det}{{\rm{Det}}}
\newcommand{\Sym}{{\rm{Sym}}}
\newcommand{\LaSa}{{\La_{S^*}}}
\newcommand{\cX}{{\cal {X}}}
\newcommand{\MHG}{{\frak {M}}_H(G)}
\newcommand{\tauMla}{\tau(M_{\lambda})}
\newcommand{\Fvur}{{F_v^{\rm {ur}}}}
\newcommand{\Lie}{{\rm {Lie}}}
\newcommand{\cL}{{\cal {L}}}
\newcommand{\cW}{{\cal {W}}}
\newcommand{\fq}{{\frak {q}}}
\newcommand{\cont}{{\rm {cont}}}
\newcommand{\SC}{{SC}}
\newcommand{\Om}{{\Omega}}
\newcommand{\dR}{{\rm {dR}}}
\newcommand{\crys}{{\rm {crys}}}
\newcommand{\hatSig}{{\hat{\Sigma}}}
\newcommand{\rdet}{{{\rm {det}}}}
\newcommand{\ord}{{{\rm {ord}}}}
\newcommand{\BdR}{{B_{\rm {dR}}}}
\newcommand{\BdRO}{{B^0_{\rm {dR}}}}
\newcommand{\Bcrys}{{B_{\rm {crys}}}}
\newcommand{\Qw}{{\mathbb{Q}}_w}
\newcommand{\barkappa}{{\bar{\kappa}}}
\newcommand{\cP}{{\Cal {P}}}
\newcommand{\cZ}{{\Cal {Z}}}
\newcommand{\oppLa}{{\Lambda^{\circ}}}
\newcommand{\Spa}{{\rm {Spa}}}
\newcommand{\tor}{{\rm {tor}}}
\newcommand{\tp}{{\rm {top}}}
\newcommand{\alg}{{\rm {alg}}}
\newcommand{\zar}{{\rm {zar}}}
\newcommand{\et}{{\rm {\acute{e}t}}}

\begin{abstract} 
  We show that the description of Deligne--Beilinson cohomology is improved by using log Hodge theory.
  We consider the log relative version of it, and also present a fundamental conjecture in log Hodge theory. 
\end{abstract}

\renewcommand{\thefootnote}{\fnsymbol{footnote}}



\section{Introduction}

\begin{para} In this paper, we show that the description of Deligne--Beilinson cohomology in Esnault--Viehweg  \cite{E} is improved by using log Hodge theory.
We also consider the log relative version of Deligne--Beilinson cohomology. 

\end{para}

\begin{para}
To treat  a smooth open algebraic variety $X$, a good method is to take a good compactification $\overline{X}$ of $X$ and put some log structure (e.g.\ log poles of differential forms) along the boundary $\overline{X} \smallsetminus X$. This is what Deligne did in Hodge II (\cite{D0}). 
  To get the mixed Hodge structure $H^m(X)=(H^m(X)_\Z, W, F)$, he defined the Hodge filtration $F$ by using a good compactification  $\overline{X}$.

However,  he did not describe the $\Z$-structure $H^m(X)_\Z=H^m(X, \Z)$  by using this compactification. In the log Hodge theory, we can use  $\overline{X}$ to describe it as  $H^m(X, \Z)\cong H^m(\overline{X}^{\log}, \Z)$. By using this description, we improve the presentation of Deligne--Beilinson cohomology in \cite{E} as is explained in Section 2.

\end{para}

\begin{para} In Section 3, we discuss how to formulate the log relative versions of Hodge II and the Deligne--Beilinson cohomology. 

\end{para}

\begin{para} {\bf Notation}. Following \cite{E}, for a scheme $X$ of finite type over $\bC$, the cohomology usually means the cohomology of the set $X(\C)$ of $\C$-points of $X$ with the classical topology. In the case where we consider the Zariski topology or the \'etale topology of the scheme $X$, we use the notation $X_{\alg}$. 
\end{para}

\section{Deligne--Beilinson cohomology}

\begin{para}\label{I1}

Let $X$ be a separated scheme of finite type over $\C$. Then for a subring $A$ of $\R$, we have Deligne--Beilinson cohomology groups $H^m_{\rm DB}(X,A(r))$ for $m, r\geq 0$, called also Deligne cohomology groups,  
 which are put in  the exact sequence
$$ \cdots \to H^m_{\rm DB}(X, A(r)) \to H^m(X, A(r)) \to H^m(X, \C)/F^r $$ $$\to H^{m+1}_{\rm DB}(X, A(r)) \to H^{m+1}(X, A(r)) \to H^{m+1}(X, \C)/F^r \to\cdots.$$
Here $A(r)$ in $H^m(X, A(r))$ denotes
 the subgroup  $A\cdot (2\pi i)^r$ of $\C$, and 
$(F^p)_p$ is the Hodge filtration on $H^m(X, \C)$. See \cite{E} Section 2.

In the case where $X$ is proper and smooth, $H^m_{\rm DB}(X, A(r))$ is defined as the cohomology of a complex of sheaves 
 as 
$$H^m_{\rm DB}(X, A(r)) := H^m(X, A(r) \to \cO_X \overset{d}\to \Omega^1_X \overset{d}\to \cdots \overset{d}\to \Omega^{r-1}_X),$$
where on the right hand side, 
$A(r)$ is put in degree $0$ in the complex, $\cO_X$ is the sheaf of holomorphic functions on $X(\C)$ and $\Omega^p_X$ is the sheaf of holomorphic $p$-forms on $X(\C)$. 

\end{para}

\begin{para}\label{I2}

Assume that $X$ is a Zariski  open set of a proper smooth scheme $\overline X$ over $\C$ in which the complement $D$ of $X$ is a divisor with normal crossings. In this paper, we give a presentation of $H^m_{\rm DB}(X, A(r))$ which is similar to the proper smooth case in \ref{I1}
by using the space ${\overline X}^{\log}$ in \cite{KN} and using the sheaves $\cO_{\overline X}^{\log}$ etc.\ on ${\overline X}^{\log}$ which appear in  log Hodge theory in \cite{KU} and \cite{KNU3}. Endow $\overline X$ with the log structure associated to $D$, and let ${\overline X}^{\log}$ be the associated space.

\end{para}

\begin{thm}\label{ast} We have the following presentation of $H^m_{\rm DB}(X, A(r))$. 
$$H^m_{\rm DB}(X, A(r))= H^m({\overline X}^{\log}, A(r) \to  \cO_{\overline X}^{\log} \overset{d}\to \Omega^{1,\log}_{\overline X}(\log D) \overset{d}\to \cdots \overset{d}\to \Omega^{r-1,\log}_{\overline X}(\log D)).$$

\end{thm}

See \cite{KN}  3.2 for the definition of $\cO_{\overline{X}}^{\log}$. The de Rham complex $\Omega^{\bullet, \log}_{\overline{X}}(\log D)= \cO_{\overline{X}}^{\log}\otimes_{\cO_{\overline{X}}} \Omega^{\bullet}_{\overline{X}}(\log D)$ on $\overline{X}$ is defined in  \cite{KN} 3.5, where it is denoted by $\omega_{\overline{X}}^{\bullet,\log}$.

\begin{para}\label{I3}
This Theorem \ref{ast}  simplifies the theory of Deligne--Beilinson cohomology. 

For example, the product structure of the Deligne--Beilinson cohomology is defined for proper smooth $X$ in  \cite{E} Section 1 in a simple way, but for $X$  in \ref{I2},  it is defined in \cite{E} Section 3 in a complicated way. By using Theorem  \ref{ast},  the product structure for $X$  in \ref{I2} is defined 
 in the same simple way as the proper smooth case in \ref{I1}.  
  See \ref{prod}. 
  Also, the Chern classes and the Chern characters  of a vector bundle  in Deligne--Beilinson cohomology explained in \cite{E} Section  8 are defined more easily by using Theorem \ref{ast} as is explained  in \ref{C1}--\ref{C3}.

We will denote the above complex $A(r) \to \cO_{\overline X}^{\log}\to \cdots \to \Omega^{r-1,\log}_{\overline X}(\log D)$ on $\overline X^{\log}$ by $A(r)_{{\rm DB}, \overline{X}}$ or simply by $A(r)_{{\rm DB}}$. 

\end{para}

\begin{para}\label{pfast} We prove Theorem \ref{ast}. 

Let $j:X\to \overline X$ and $j^{\log}:X\to {\overline X}^{\log}$ be the inclusion maps and let  $\theta: {\overline X}^{\log}\to {\overline X}$  be the canonical map. 
Then $j=\theta\circ j^{\log}$.

 By definition, $H^m_{\rm DB}(X, A(r))=H^m(\overline X, C)$,  where $C$ is the mapping fiber ($-1$ translation of the mapping cone) of
 $$\Omega^{\geq r}_{\overline X}(\log D)\oplus Rj_*A(r) \to Rj_*\Omega^{\bullet}_X.$$ Here $\Omega^{\geq r}_{\overline X}(\log D)$ denotes the part of degree $\geq r$ of 
 $\Omega^{\bullet}_{\overline X}(\log D)$. On the other hand,  $A(r)_{{\rm DB}}$  is quasi-isomorphic to the mapping fiber of $\Omega_{\overline X}^{\geq  r,\log}(\log D)\oplus A(r) \to \Omega_{\overline X}^{\bullet, \log}(\log D)$. 
 Hence it is sufficient to prove that the vertical arrows in the following commutative diagram are isomorphisms.
 $$\begin{matrix}  \Omega_{\overline X}^{\geq r}(\log D) & \to & Rj_*\Omega^{\bullet}_X & \leftarrow & Rj_*A(r)\\
 \downarrow && \uparrow && \uparrow \\
 R\theta_* \Omega_{\overline X}^{\geq r,\log}(\log D) & \to& R\theta_*\Omega^{\bullet,\log}_{\overline X}(\log D) &\leftarrow& R\theta_*A(r)\end{matrix}$$
 Here the middle and the right vertical arrows are defined by $j=\theta\circ j^{\log}$. 
 The left and the middle vertical arrows are isomorphisms by the following (1) which we apply to $\cF=\Omega^p_{\overline X}(\log D)$ together with the well-known isomorphism 
 $\Omega^{\bullet}_{\overline X}(\log D) \cong Rj_*\Omega^{\bullet}_{X}$ for the middle. 
 The right vertical arrow is an isomorphism by  the following (2) which we apply to the inverse image $\cF$ of $A(r)$ on ${\overline X}^{\log}$. 
 
(1) (\cite{KU} Proposition 2.2.10.) For a vector bundle $\cF$ on $\overline X$, we have $\cF \overset{\cong}\to R\theta_*(\cO_{\overline X}^{\log}\otimes_{\cO_X} \cF)$. 
  
(2) (\cite{KN} Remark 1.5.1.) For a locally constant sheaf of abelian groups $\cF$ on ${\overline X}^{\log}$, we have $\cF\overset{\cong}\to Rj^{\log}_*j^{\log*}\cF$.

\end{para}

\begin{para}\label{prod} For $p,p'\geq 0$, we define a multiplication $$\cup: A(p)_{\rm DB} \otimes A(p')_{\rm DB}\to A(p+p')_{\rm DB}$$ as follows, in the same way as in \cite{E} Section 1: $x\cup y$ is $xy$ if $x$ is of degree $0$, $x\wedge dy$ if $x$ is of degree $>0$ and $y$ is of degree $p'$, and is $0$ otherwise. 
The multiplication 
$$\cup: H^m_{\rm DB}(X, A(p)) \otimes H^n_{\rm DB}(X, A(p')) \to H^{m+n}_{\rm DB}(X, A(p+p'))$$ is induced by this.

\end{para}

\begin{para} 
Let $X$ be a separated scheme of finite type over $\C$.
We have  a presentation of $H^m_{\rm {DB}}(X, A(r))$ using Theorem \ref{ast} for $X$ in \ref{I2} and using  a hyper-covering. 

As in (8.2.1) of \cite{D}, 
there is a  diagram
$$\begin{matrix}   Y _{\bullet} &\overset{\subset}\to &  \overline Y_{\bullet}\\
\downarrow &&\\
X&&\end{matrix}$$
in which 
$Y_{\bullet}$ is the complement of a normal crossing divisor $D_{\bullet}$ in the proper simplicial smooth  scheme $\overline Y_{\bullet}$ and $Y_{\bullet}\to X$ is a proper hyper-covering. Using this diagram, Deligne gives the mixed Hodge structure on $H^m(X, \Z)$ which is isomorphic to $H^m(Y_{\bullet},  \Z)$.  We have the presentation of  Deligne--Beilinson cohomology by log Hodge theory as $$H^m_{\rm DB}(X, A(r))= H^m({\overline Y}^{\log}_{\bullet}, A(r)_{{\rm DB}, \overline Y_{\bullet}}).$$

\end{para}

\begin{para}\label{C1} In \ref{C1}--\ref{C3}, we show that the Chern classes and the Chern characters of vector bundles in Deligne--Beilinson cohomology are explained in nice ways by our log method.

For $X$ in \ref{I2}, the first Chern class mapping for  line bundles 
on $X_{\alg}$ is defined as

\smallskip

(1) $\mathrm{Pic}\,(X_{\alg})=H^1(X_{\alg,\et}, \cO_{X,\alg}^\times)\cong H^1(\overline X_{\alg,\et}, M_{\overline X,\alg}^{\gp}) \to H^2_{\rm DB}(X, \Z(1)).$

\smallskip

Here the  isomorphism in the middle is obtained from $j_{\alg, \et,*}(\cO_{X_{\alg}}^{\times})= M_{X_{\alg}}^{\gp}$ and $R^1j_{\alg,\et,*}(\cO_{X_{\alg}}^{\times})=0$. 
The last map in (1) is induced by $\log: \theta^{-1}(M_{\overline X}^{\gp})\to \cO^{\log}_{\overline X}/\Z(1)$. 

\end{para}

\begin{para}\label{C2} For a separated scheme $X$ of finite type over $\C$, the first Chern class mapping for  line bundles
on $X_{\alg}$  is defined as 
$${\mathrm {Pic}}\,(X_{\alg}) \to {\mathrm {Pic}}\,(Y_{\alg, \bullet})\cong H^1({\overline Y}_{\bullet, \alg, \et}, M^{\gp}_{Y_{\bullet, \alg}})$$ 
$$\to H^1({\overline Y}^{\log}_{\bullet}, \cO_{{\overline Y}_{\bullet}}^{\log}/\Z(1))= H^2({\overline Y}^{\log}_{\bullet}, \Z(1)\to \cO_{{\overline Y}_{\bullet}}^{\log})=H^2_{\rm DB}(X, \Z(1))$$
by using the hyper-covering. 
\end{para}

\begin{para}\label{C3} For a vector bundle $\cF$ on $X_{\alg}$, the  Chern class and the Chern character of  $\cF$ in the Deligne--Beilinson cohomology are defined by using the above first Chern classes of  line bundles in the standard way.

\end{para}

\section{Log relative versions}\label{s:logver} 

In \ref{D1}--\ref{status} (resp.\ \ref{DB1}--\ref{Mlog}), we consider how to formulate the logarithmic relative version of Hodge II (\cite{D0}) (resp.\  of the theory of Deligne--Beilinson cohomology).

\begin{para}\label{D00} In this Section \ref{s:logver}, let $f: X\to S$ be a projective and log smooth saturated morphism between fs log complex analytic spaces. 

The present $X\to S$ generalizes $\overline{X}\to \Spec(\C)$ of Section 2. 

\end{para}

\begin{para}\label{D0} Let $m\geq 0$. We have:

(1) Concerning the induced map $f^{\log}: X^{\log}\to S^{\log}$, $R^mf^{\log}_*\Z$ is a locally constant sheaf of finitely generated $\Z$-modules on $S^{\log}$ (\cite{KajiwaraNakayama}).

(2) We have a canonical isomorphism $$\cO_S^{\log} \otimes_{\Z} R^mf^{\log}_*\Z\cong \cO_S^{\log}\otimes_{\cO_S} H^m_{\dR}(X/S)$$
(\cite{IKN} Theorem (6.3)), where $H^m_{\dR}(X/S)= R^mf_*(\omega^{\bullet}_{X/S})$ with $\omega^{\bullet}_{X/S}$ being the logarithmic de Rham complex of $X$ relative to $S$. 
By this, the $\cO_S$-module $H^m_{\dR}(X/S)$ is locally free of finite rank.
\end{para}
\begin{para}\label{D1} In \ref{D1}--\ref{status}, 
let $f: X\to S$ be as in \ref{D00} and we assume 
 that 
the following condition \ref{D1}.1  is satisfied.

\medskip
  
\noindent {\bf \ref{D1}.1.}   Locally on $X$ and on $S$, $X$ is isomorphic over $S$ to an open subspace of $Y\times \C^r$, where $Y$ is an fs log analytic space over $S$ which is log smooth and vertical over $S$ and $\C^r$ is endowed with the log structure generated by the coordinate functions. 

\medskip

 This condition is assumed to have a reasonable weight filtration. This condition is automatically satisfied if $f$ is vertical.

  \end{para}

\begin{para} \label{pLMH} We consider the $4$-ple $$H^m(X/S):=(H^m(X/S)_\Z, H^m(X/S)_\cO, W, F)$$ 
defined as follows. Let  $H^m(X/S)_\Z$ be the locally constant sheaf on $S^{\log}$ given by 
$$H^m(X/S)_\Z:= R^mf^{\log}_*\Z.$$ Let $H^m(X/S)_{\cO}$ be the vector bundle $H^m_{\dR}(X/S)$ on $S$. We have a canonical isomorphism $\cO_S^{\log}\otimes_{\Z} H^m(X/S)_\Z\cong \cO_S^{\log}\otimes_{\cO_S} H^m(X/S)_{\cO}$ (\ref{D0}).

Define the  Hodge filtration $F$  on $H^m(X/S)_{\cO}$ as follows: Let $F^p$ be the image of $R^mf_*(\omega^{\geq p}_{X/S})$.

Define the weight filtration $W$ on $H^m(X/S)_\Q:=\Q\otimes_{\Z} H^m(X/S)_\Z$ as follows:
Let $X'$ be the fs log analytic space whose underlying analytic space is $X$ and whose log structure is the vertical part $\{a\in M_X\;|\; a|b \;\text{for some}\;b\in M_S\}$ of the log structure of $X$ relative to $S$, and let 
 $\theta:X^{\log}\to (X')^{\log}$ and $\eta:(X')^{\log}\to S^{\log}$ 
be the 
canonical maps. Let $W_w$ be the image of $R^m\eta_* \tau_{\leq w-m} R\theta_*\Q$, where $\tau_{\leq \bullet}$ is the canonical truncation functor. 

\end{para}

\begin{para}\label{D2} By the definition of the weight filtration, $H^m(X/S)$ is of weight $\geq m$. If $X/S$ is vertical, it is of pure weight $m$ (in this case, $X=X'$). 

\end{para}

The following is one of the main problems in log Hodge theory. 

\begin{conj}
\label{conj} $(1)$ The Hodge to de Rham spectral sequence $$E^{p,q}_1= R^qf_*(\omega^p_{X/S}) \Rightarrow H^{p+q}_{\dR}(X/S)$$ degenerates at $E_1$ and all $R^qf_*(\omega^p_{X/S})$ are locally free of finite rank as $\cO_S$-modules. 
Hence the $4$-ple $H^m(X/S)$ is a pre-log mixed Hodge structure (pre-LMH) ({\rm \cite{KU} 2.6}, {\rm \cite{KNU3} 1.3.1}) on $S$. 

$(2)$ The pre-LMH $H^m(X/S)$ is a log mixed Hodge structure (LMH) ({\rm \cite{KU} 2.6}, {\rm \cite{KNU3} 1.3.1}) on $S$. 
\end{conj}

\begin{para}
\label{status} Concerning Conjecture \ref{conj}:

  By \cite{SZ}, it holds if $X$ has a semistable reduction over a log pointed disk $S$.

  By \cite{KMN}, it holds if $f$ is vertical and $S$ is log smooth. 
  
   By \cite{FN}, it holds if $f$ is vertical and the log rank of $S \le 1$.

  The part (1) of \ref{conj} holds
  if every  fiber of $f$ is algebraic. Here, the underlying analytic space of every fiber of $f$ is algebraic because $f$ is projective, but this assumption means that the log structure of the fiber is also algebraic; the authors do not know whether it is automatic or not. 
  In fact, the statement is reduced to that for the fibers by the usual argument by the lengths (see \cite{IKNe}), and if the fiber is algebraic, then the statement for the fiber is reduced to \cite{IKN} Theorem (7.1), 
in which $X$ and $S$ are assumed to be algebraic. 

  We plan to discuss how to attack the conjecture \ref{conj} in the sequel of our papers \cite{KNU3} and \cite{KNU4}. 
\end{para}

\begin{para}\label{DB1} In the rest of this paper, let $f: X\to S$ be as in \ref{D00}.  Let  $\omega^{p, \log}_{X/S}$ be the sheaf $ \cO_X^{\log} \otimes_{\cO_X} \omega^p_{X/S}$ on $X^{\log}$. 
Let $a: S^{\log}\to S$ be the canonical map and let $$g:= a \circ f^{\log}: X^{\log}\overset{f^{\log}}\to S^{\log}\overset{a}\to S.$$

We define the log DB cohomology of $X/S$ by $$H^m_{\log\mathrm{DB}}(X/S, A(r)):=R^mg_*(A(r)_{{\rm DB}, X/S}), $$ 
where $A(r)_{{\rm DB}, X/S}$ denotes the complex $A(r) \to \cO_X^{\log} \overset{d}\to \omega_{X/S}^{1, \log}\overset{d}\to \cdots \overset{d}\to \omega_{X/S}^{r-1, \log}$. 

\end{para}

\begin{para}\label{fXS2}
 In Proposition \ref{Ext} below, we will assume that $X/S$ is vertical and that the log rank of $S$ is $\leq 1$. The following are such examples. 
\medskip

{\it Examples.} Let $C$ be a smooth algebraic curve  over $\C$ and let $\frak X\to C$ be projective, generically smooth, and of semi-stable reduction. Let $\Sig$ be a  finite closed subset   of $C$ outside which $\frak X$ is smooth, endow $C$ with the log structure associated to $\Sig$, and endow $\frak X$ with the log structure associated to the inverse image of $\Sig$ which is a divisor with normal crossings.

(i) Let $X=\frak X$, $S=C$. Then 
$f:X\to S$ satisfies the assumption of this \ref{fXS2}.

(ii) Let  $s=S$ be a point in $\Sig$, and let $X$ be the fiber of $s$ in $\frak X$. Then with the induced  log structures, $f: X\to S$ satisfies the assumption of \ref{fXS2}.

  Note that Conjecture \ref{conj} is already solved in these cases (\ref{status}). 
  
  \end{para}

\begin{prop}\label{Ext}  Assume that $X/S$ is vertical. 
Let $m, r\geq 0$. Let 
$H^m_{\log\mathrm{B}}(X/S, \Z(r)):= R^mg_*(\Z(r))$. 
  Assume that $S$ is of log rank $\leq 1$. Assume $m+1-2r\leq 0$. 

$(1)$ We have a commutative diagram with exact rows 
$$\begin{matrix} 0 & \to &  \cE xt_{\cH}(\Z, H^m(X/S)(r)) & \to & H^{m+1}_{\log\mathrm{DB}}(X/S,\Z(r))& \to & \cH om_{\cH}(\Z, H^{m+1}(X/S)(r))&\to& 0\\
&&\downarrow  &&\downarrow&& \cap &&\\
0&\to &\cE xt_{\cB}(\Z, H^m(X/S)(r)_\Z) & \to & H^{m+1}_{\log\mathrm{B}}(X/S, \Z(r))& \to & \cH om_{\cB}(\Z,  H^{m+1}(X/S)(r)_\Z)&\to & 0.\end{matrix}$$
Here $\cH om_{\cH}(\dots)$ (resp.\ $\cH om_{\cB}(\dots)$) denotes the sheaf $U\mapsto \Hom_{\cH(U)}(\dots)$ (resp.\ $U\mapsto \Hom_{\cB(a^{-1}(U))}(\dots)$ on $S$ with $\cH(U)$ (resp.\ $\cB(a^{-1}(U))$) being the category of
 pre-LMH on $U$ (resp.\  sheaves of abelian groups on $a^{-1}(U)\subset S^{\log}$) for an open set $U$ of $S$, and $\cE xt_{\cH}(\dots)$ (resp.\ $\cE xt_{\cB}(\dots)$) denotes the sheaf on $S$ associated to the pre-sheaf $U\mapsto \Ext_{\cH(U)}(\dots)$ (resp.\ $U\mapsto \Ext_{\cB(a^{-1}(U))}(\dots)$).

$(2)$ The left vertical arrow is surjective. 
The kernel of the left vertical arrow and the kernel of the middle vertical arrow are identified with $(a_*(H^m(X/S)(r)_\Z)) \bs H^m_{\mathrm{dR}}(X/S)/F^rH^m_{\mathrm{dR}}(X/S)$.

$(3)$ If  $m+1-2r<0$,  $\cH om_{\cH}(\Z, H^{m+1}(X/S)(r))$ is finite  and hence 
 $$\Q\otimes \cE xt_{\cH}(\Z, H^{m}(X/S)(r)) \overset{\cong}\to H^{m+1}_{\log\mathrm{DB}}(X/S,\Q(r)).$$

$(4)$ If $s\in S$ and if $t\in S^{\log}$ lies over $s$, the stalk of $\Q\otimes \cE xt_{\cB}(\Z, H^m(X/S)(r)_\Z)$ 
 at $s$ is isomorphic to $H^m(X/S)(r)_{\Q,t}/N_tH^m(X/S)(r)_{\Q,t}$ (resp.\ $0$) if the log structure is of rank $1$ (resp.\ $0$) at $s$. Here $N_t$ is the monodromy logarithm at $t$. 

\end{prop}

To prove Proposition \ref{Ext}, we use

\begin{lem}\label{Ext2} Let $H=(H_\Z, W, F)$ be a pre-LMH on an fs log analytic space $S$. Let $C$ be the complex $H_\Z \to \cO^{\log}_S \otimes_{\cO_S} H_\cO/F^0$ of sheaves on $S$, where $H_\Z$ is put in degree $0$. 

$(1)$ If $H$ is of weights $\leq 0$, we have $\cH om_{\cH}(\Z, H)\cong R^0a_*C$. 

$(2)$ If  $H$ is of weights $\leq -1$, 
$\cE xt_{\cH}(\Z, H)\cong R^1a_*C$.

\end{lem}

\begin{pf} (1) follows from $R^0a_*C= \text{Ker}(a_*H_\Z\to H_{\cO}/F^0)$. 

 We prove (2). Note that $R^1a_*C$ is the sheaf of isomorphism classes of pairs $(L, \la)$, where $L$ is a sheaf of abelian groups on $S$ with an exact sequence $0\to H_\Z \overset{v}\to L \to \Z \to 0$ and $\la$ is a homomorphism $L\to \cO_S^{\log}\otimes_{\cO_S} H_{\cO}/F^0$
 such that $\la \circ v$ is the canonical map. For such an $(L, \la)$, we have an extension $(L, \tilde W, \tilde F)$, where $\tilde W$ is the evident one and $\tilde F$ is the extension of $F$ to $\cO_S \otimes L$ characterized by the property that $\cO_S^{\log}\otimes_{\cO} \tilde F^0$ is generated by $F^0$ of $H$ and $x+ \la(x)$ ($x\in L$). 
\end{pf}

\begin{para}\label{pfExt} We prove Proposition \ref{Ext}. The lower exact sequence in (1) comes from the spectral sequence
$$E^{p,q}_2  = R^pa_*(H^q(X/S)(r)_\Z) \Rightarrow H^{p+q}_{\log\mathrm{B}}(X/S, \Z(r)).$$
The upper exact sequence in (1) comes from the spectral sequence
$$E^{p,q}_2  = R^pa_*[H^q(X/S)(r)_\Z \to \cO_S^{\log} \otimes_{\cO_S} H^q_{\mathrm{dR}}(X/S)/F^r]\Rightarrow H^{p+q}_{\log\mathrm{DB}}(X/S, \Z(r)),$$ 
as follows. By Lemma \ref{Ext2}, $E^{0,m+1}_2=\cH om_{\cH}(\Z, H^{m+1}(X/S)(r))$ and $E^{1,m}_2= \cE xt_{\cH}(\Z, H^m(X/S)(r))$. 
Furthermore, we have 

\medskip

\ref{pfExt}.1  $R^ja_*(\cO_S^{\log} \otimes_{\cO_S} H^q_{\mathrm{dR}}(X/S)/F^r)=0$ for $j>0$ 

\medskip
\noindent
by \ref{pfast} (1),  and we have $R^ja_*=0$ for $j>1$ because $S$ is of log rank $\leq 1$, and hence $E^{p,q}_2=0$ for $p\geq 2$. Hence we have the upper exact sequence in (1). 

(2) follows from \ref{Ext2} (2) and \ref{pfExt}.1. 

(3) is a consequence of the exactness of the upper row in (1).

(4) is by generality of local systems.

\end{para}

\begin{para} There are differences between the categories $\cH(U)$ of pre-LMH and the categories of LMH (\cite{KU} 2.6, \cite{KNU3} 1.3.1). An LMH is a pre-LMH which satisfies the admissibility of local monodromy,  the Griffiths transversality of the pullbacks to log points, and the condition that sufficiently twisted specializations give mixed Hodge structures in the usual sense  (these are conditions (1), (2), (3) in \cite{KNU3} 1.3.1, respectively).

Hence $\cE xt$ for LMH can be smaller than the $\cE xt_{\cH}$ in the above proposition.
\end{para}

\begin{para} Let the assumption be as in Proposition \ref{Ext}. 
In \cite{KNU3} Sections 5 and 6, we considered the moduli space of extensions $0\to H^m(X/S)(r)\to E\to \Z\to 0$ of LMH. By Proposition \ref{Ext}, this space  is regarded as a part of  
$H^{m+1}_{\log\mathrm{DB}}(X/S, \Z(r))$.

\end{para}

\begin{para} Let $\frak X \to C\supset \Sig$ be as in \ref{fXS2} and assume that $C$ is proper.  For the inverse image $U$ of $C\smallsetminus \Sig$ in $\frak X$,  an element $k$ of the $K$-group $K_{2r-1-m}(U)$ with $m< 2r-1$ gives a  log mixed Hodge structure $E_k$ 
on $C$ with an exact sequence $0\to H^m(\frak X/C)(r)
\to E_k \to \Z\to 0$ (cf. \cite{KNU4} Section 7.2). By Proposition \ref{Ext}, how $E_k$ degenerates at $s\in \Sig$ is related to the image of $k$ under the composition   $$K_{2r-1-m}(U) \to H^{m+1}_{\mathrm{DB}}(U, \Z(r))\cong H^{m+1}(\frak X^{\log}, \Z(r)_{\mathrm DB})\to H^{m+1}_{\log\mathrm{DB}}(\frak X_s/s, \Z(r)),$$
where the first arrow is the Beilinson regulator and  $\frak X_s$ 
denotes the fiber of $s$ in $\frak X$.  

We expect that log Deligne--Beilinson cohomology is useful for understanding of 
how elements of $K$-groups and algebraic cycles behave in degeneration. 
\end{para}

\begin{para}\label{Mlog} 

We expect that in a suitable algebraic situation over $\C$, we have a cycle map (the Hodge realization map) from the log motivic cohomology  in \cite{P}, \cite{B} to the log Deligne--Beilinson cohomology. 
\end{para}

\end{document}